\documentclass[article]{elsarticle}

\usepackage{hyperref}
\usepackage{amsfonts,amssymb}
\usepackage{mathtools}
\usepackage{multirow}%

\journal{Journal of \LaTeX\ Templates}

\newtheorem{thm}{Theorem}

\newdefinition{rem}{Remark}

\begin{document}

\begin{frontmatter}
\title{Optimal stopping via reinforced  regression} 
\tnotetext[mytitlenote]{This work was supported by   RSF grant 19-71-30020}
\author[address1,address2]{D. Belomestny\corref{mycorrespondingauthor}}
\ead[url]{www.uni-due.de/~hm0124}
\ead{denis.belomestny@uni-due.de}
\author[address3]{J. Schoenmakers}
\author[address2,address3]{V. Spokoiny}
\author[address2]{B. Zharkynbay}
\cortext[mycorrespondingauthor]{Corresponding author}
\address[address1]{Duisburg-Essen University, Essen}
\address[address2]{National University Higher School of Economics, Moscow}
\address[address3]{WIAS, Berlin}
\begin{abstract}
In this note we propose a new approach towards solving numerically optimal
stopping problems via reinforced regression based Monte Carlo algorithms. The
main idea of the method is to reinforce standard linear regression algorithms in
each backward induction step by adding new basis functions based on previously
estimated continuation values. The proposed methodology is illustrated by
several numerical examples from mathematical finance.
\end{abstract}
\end{frontmatter}


\section{Introduction}

A discrete time optimal stopping problem can be efficiently solved in low
dimensions, for instance by tree methods or by using deterministic numerical
methods for the corresponding partial differential equation. However, many
optimal stopping problems arising in applications (see e.g. \cite{Gl}) involve
high dimensional underlying processes and this made it necessary to develop
Monte Carlo methods for solving such problems. Solving optimal stopping
problems via Monte Carlo is a challenging task, because this typically
requires backward dynamic programming that for long time was thought to be
incompatible with forward structure of Monte Carlo methods. In recent years
much research was focused on the development of efficient methods to compute
approximations to the value functions or optimal exercise policy. Eminent
examples include the functional optimization approach of \cite{A}, the mesh
method of \cite{BG}, the regression-based approaches of \cite{Car},
\cite{J_LS2001}, \cite{J_TV2001}, \cite{egloff2005monte} and \cite{B1}. The
most popular type of algorithms are with no doubt the regression ones. In
fact, in many practical pricing problems, the low-degree polynomials are
typically used for regression (see \cite{Gl}). The resulting least squares
problem has a relatively small number of unknown parameters. However, this
approach has an important disadvantage - it may exhibit too little flexibility
for modelling highly non-linear behaviour of the exercise boundary.
Higher-degree polynomials can be used, but they may contain too many
parameters and, therefore, either over-fit the Monte Carlo sample or prohibit
parameter estimation because the number of parameters is too large. In this
note a regression based Monte Carlo approach is developed for building sparse
regression models at each backward step of the dynamic programming algorithm.
This enables estimating the value function with virtually the same cost as the
standard regression algorithms based on low degree polynomials but with higher
precision. The additional basis functions are constructed specifically for the
optimal stopping problem at hand without using a fixed predefined finite
dictionary. Specifically, the new basis functions are learned during the
backward induction via incorporating information from the preceding backward
induction step. Our algorithm may be viewed as a method of constructing sparse
nonlinear approximations (in terms of their dependence on Monte Carlo paths)
of the underlying value function and in this sense it extends the literature
on nonlinear learning type algorithms for optimal stopping problems, see, for
example, the recent paper \cite{becker2018deep} and references therein.

The structure of the paper is as follows. After recalling basic facts on
American options and settling the main setup in Section~\ref{mainsetup}, the
reinforced procedure is presented in Section~\ref{boost}. The numerical
performance is studied in Section~\ref{num}.

\section{Main setup}

\label{mainsetup} A general class of optimal stopping problems respectively,
can be formulated with respect to an underlying $\mathbb{R}^{d}$-valued Markov
process $(X_{t},\,0\leq t\leq T)$ defined on a filtered probability space
$(\Omega,\mathcal{F},(\mathcal{F}_{t})_{0\leq t\leq T},\mathrm{P})$. The
process $(X_{t})$ is assumed to be adapted to a filtration $(\mathcal{F}%
_{t})_{0\leq t\leq T}$ in the sense that each $X_{t}$ is $\mathcal{F}_{t}$
measurable. Recall that each $\mathcal{F}_{t}$ is a $\sigma$ -algebra of
subsets of $\Omega$ such that $\mathcal{F}_{s}\subseteq\mathcal{F}%
_{t}\subseteq:\mathcal{F}$ for $s\leq t.$ Henceforth we restrict our selves to
the case where only a finite number of stopping opportunities $0<t_{1}%
<t_{2}<\ldots<t_{\mathcal{J}}=T$ are allowed. We now consider~the
pre-specified reward process $g_{j}(Z_{j})$ in terms of the Markov chain
\[
Z_{j}:=X_{t_{j}},\quad j=1,\ldots,\mathcal{J},
\]
for some given functions $g_{1},\ldots,g_{\mathcal{J}}$ mapping $\mathbb{R}%
^{d}$ into $[0,\infty).$ Let $\mathcal{T}_{j}$ denote the set of stopping
times taking values in $\{j,j+1,\ldots,\mathcal{J}\}$ and consider the optimal
stopping problems of the form
\begin{equation}
V_{j}(x)=\sup_{\tau\in\mathcal{T}_{j}}\mathsf{E}[g_{\tau}(Z_{\tau}%
)|Z_{j}=x],\quad x\in\mathbb{R}^{d},\label{stop}%
\end{equation}
In (\ref{stop}) we have to read $\mathcal{T}_{0}:=\mathcal{T}_{1}$ for $j=0.$
A common feature of many approximation algorithms for optimal stopping
problems is that they deliver estimates $C_{N,1}(x),\ldots,C_{N,\mathcal{J}%
-1}(x)$ for the so-called continuation functions:
\begin{equation}
C_{j}(x):=\mathsf{E}[V_{j+1}(Z_{j+1})|Z_{j}=x],\quad j=1,\ldots,\mathcal{J}%
-1.\label{dyn0}%
\end{equation}
Here the index $N$ indicates that the above estimates are based on a set of
$N$ independent \textquotedblleft training\textquotedblright\ trajectories
\begin{equation}
(Z_{1}^{(i)},\ldots,Z_{\mathcal{J}}^{(i)}),\text{ \ \ }i=1,\ldots
,N,\label{train}%
\end{equation}
all starting from one point. In the case of the so-called regression methods,
the estimates for (\ref{stop}) and (\ref{dyn0}) are obtained via the
application of \textit{Dynamic Programming Principle}:
\begin{align}
C_{j}(x)  &  =\mathsf{E}[V_{j+1}(Z_{j+1})|Z_{j}=x], \quad V_{j}(x)
=\max\left(  g_{j}(x),C_{j}(x)\right)  ,\text{ \ \ }1\leq j\leq\mathcal{J}%
-1,\nonumber
\end{align}
with $V_{\mathcal{J}}(x) =g_{\mathcal{J}}(x),$ $C_{\mathcal{J}}(x)=0,$
combined with nonparametric regression.

In the setting of Tsitsiklis-van Roy \cite{J_TV2001}, this regression
algorithm can be described as follows. First initialize $C_{N,\mathcal{J}%
}(x)\equiv0.$ Suppose that for some $1\leq j<\mathcal{J},$ an estimate
$C_{N,j+1}(x)$ for $C_{j+1}(x)$ is already constructed. Then in the $j$th step
one needs to estimate the conditional expectation
\begin{equation}
\mathsf{E}[V_{N,j+1}(Z_{j+1}))|Z_{j}=x],\label{regr_aim}%
\end{equation}
where $V_{N,j+1}(x)$ $=$ $\max\left(  g_{j+1}(x),C_{N,j+1}(x)\right) .$ This
can be done by performing nonparametric regression (linear or nonlinear) on
the set of paths
\begin{equation}
(Z_{j}^{(i)},V_{N,j+1}(Z_{j+1}^{(i)})),\quad i=1,\ldots,N,\label{regr_aim1}%
\end{equation}
due to a family of basis functions resulting in the estimate $C_{N,j}(x).$

In the method of Longstaff-Schwartz \cite{J_LS2001}, one constructs the
estimates, $C_{N,j}^{\text{LS}}$ say, by regression using an interleaving set
of ``dummy cash-flows'' $\widehat{V}_{j}^{(i)}$ in the following way. First
initialize, besides $C_{\mathcal{J}}^{\text{LS}}\equiv0$, $\widehat
{V}_{\mathcal{J}}^{(i)}:=g_{\mathcal{J}}(Z_{\mathcal{J}}^{(i)}),$ $i=1,...,N.$
Once $C_{N,j+1}^{\text{LS}}$ and $\widehat{V}_{j+1}^{(i)}$ are constructed for
$j+1\leq\mathcal{J},$ compute the regression estimate $C_{N,j}^{\text{LS}}$
with respect to some set of basis functions via (\ref{regr_aim1}) with
$V_{N,j+1}(Z_{j+1}^{(i)})$ replaced by $\widehat{V}_{j+1}^{(i)}.$ Next update
\[
\widehat{V}_{j}^{(i)}=\left\{
\begin{tabular}
[c]{l}%
$g_{j}(Z_{j}^{(i)}),$ \ \ $g_{j}(Z_{j}^{(i)})\geq C_{N,j}^{\text{LS}}%
(Z_{j}^{(i)});$\\
$\widehat{V}_{N,j+1}^{(i)},$ \ \ $g_{j}(Z_{j}^{(i)})<C_{N,j}^{\text{LS}}%
(Z_{j}^{(i)}),$%
\end{tabular}
\right.
\]
for $i=1,...,N$ (see also \cite{Gl}).


Given the estimates $C_{N,1}(x),\ldots,C_{N,\mathcal{J}-1}(x)$ (Tsitsiklis-van
Roy or Longstaff-Schwartz), we next may construct a lower bound (low biased
estimate) for $V_{0}$ using the (generally suboptimal) stopping rule:
\[
\tau_{N}=\min\bigl\{1\leq j\leq\mathcal{J}:g_{j}(Z_{j})\geq C_{N,j}%
(Z_{j})\bigr\},
\]
with $C_{N,\mathcal{J}}\equiv0$ by definition. Indeed, fix a natural number
$N_{\text{test}}$ and simulate $N_{\text{test}}$ new independent trajectories
of the process $Z.$ A low-biased estimate for $V_{0}$ can be then constructed
as
\begin{equation}
V_{0}^{N_{\text{test}},N}=\frac{1}{N_{\text{test}}}\sum_{r=1}^{N_{\text{test}%
}}g_{\tau_{N}^{(r)}}\bigl(Z_{\tau_{k}^{(r)}}^{(r)}\bigr)\label{form1}%
\end{equation}
with
\begin{equation}
\tau_{N}^{(r)}=\min\Bigl\{1\leq j\leq\mathcal{J}:g_{j}(Z_{j}^{(r)})\geq
C_{N,j}(Z_{j}^{(r)})\Bigr\}.\label{form2}%
\end{equation}

\section{Reinforced regression algorithms}

\label{boost}

In this section we outline our methodology for estimating the solution to
(\ref{stop}) at time $t=0,$ based on a set of training trajectories
(\ref{train}). In this respect, as a novel ingredient, we will reinforce the
standard regression procedures by learning and incorporating new basis
functions on the backward fly. As a canonical example one may consider
incorporation of $V_{N,j} $ as a basis function in the regression step of
estimating $C_{j-1}.$ Other possibilities are, for example, certain (spatial)
derivatives of $V_{j},$ or functions directly related to the underlying
exercise boundary at time $j,$ for example $1_{\left\{  g_{j}-C_{N,j}\right\}
}.$ In general one may choose a (typically small) number of suitable
reinforcing basis functions at each step.

\subsection{Backward reinforcement of regression basis}

Let us suppose that we have at hand some fixed and a computationally cheep
system of basis functions $\left(  \psi_{1}(x),\ldots,\psi_{K}(x)\right)  .$
We now extend this basis at each backward regression step $j-1$ with an
additional and sparse set of new functions $v_{1}^{N,j-1},\ldots,v_{b}%
^{N,j-1}$ that are constructed in the preceding backward step $j, $ on the
given training paths. The main idea is that the so constructed basis delivers
more accurate regression estimate $C_{N,j-1}$ of the continuation function
$C_{j-1},$ compared to the original basis, and at the same time remains cheap.

\subsection{Backward reinforced regression algorithm}

\label{psalg} Based on the training sample (\ref{train}), we propose a
reinforced backward algorithm that in pseudo-algorithmic terms works as
follows. At time $\mathcal{J}$ we initialize as $C_{N,\mathcal{J}}(x)=0.$
Suppose that for $j<\mathcal{J},$ $C_{N,j}$ is already constructed in the
form
\[
C_{N,j}(x)=\sum_{k=1}^{K}\gamma_{k}^{N,j}\psi_{k}(x)+\sum_{k=1}^{b}%
\gamma_{k+K}^{N,j}\nu_{k}^{N,j}(x)\text{ \ \ for some \ }\gamma^{N,j}%
\in\mathbb{R}^{K+b}.
\]
For going from $j>0$ down to $j-1,$ define the new reinforced regression basis
via
\begin{equation}
\Psi^{N,j-1}(x):=\left(  \psi_{1}(x),\ldots,\psi_{K}(x),\nu_{1}^{N,j-1}%
(x),\ldots,\nu_{b}^{N,j-1}(x)\right) \label{psij}%
\end{equation}
(as a row vector) due to a choice of the set of functions $(\nu_{1}%
^{N,j-1},\ldots,\nu_{b}^{N,j-1})$ based on the previously estimated
continuation value $C_{N,j}$. For example, we might take $b=1$ and consider
the function
\begin{equation}
\nu_{1}^{N,j-1}(x)=\max(g_{j}(x),C_{N,j}(x)).\label{vb}%
\end{equation}
Then consider the $N\times\left(  K+b\right)  $ design matrix $\mathcal{M}%
^{j-1}$ with entries.
\begin{equation}
\mathcal{M}_{mk}^{j-1}:=\Psi_{k}^{N,j-1}(Z_{j-1}^{(m)}),\text{ \ \ }%
m=1,\ldots,N,\text{ }k=1,\ldots,K+b,\label{mpsi}%
\end{equation}
and the (column) vector%
\begin{align}
\mathcal{V}_{j}  &  =\left(  V_{N,j}(Z_{j}^{(1)}),\ldots,V_{N,j}(Z_{j}%
^{(N)})\right)  ^{\top}\label{vp1}\\
&  =\left(  \max(g_{j}(Z_{j}^{(1)}),C_{N,j}(Z_{j}^{(1)})),\ldots,\max
(g_{j}(Z_{j}^{(N)}),C_{N,j}(Z_{j}^{(N)}))\right)  ^{\top}.\nonumber
\end{align}
Next compute and store%
\begin{equation}
\gamma^{N,j-1}:=\left(  \left(  \mathcal{M}^{j-1}\right)  ^{\top}%
\mathcal{M}^{j-1}\right)  ^{-1}\left(  \mathcal{M}^{j-1}\right)  ^{\top
}\mathcal{V}_{j},\label{gjm1}%
\end{equation}
and then set%
\begin{align}
C_{N,j-1}(x)  &  =\Psi^{N,j-1}(x)\gamma^{N,j-1}\label{cjm1}\\
&  =\sum_{k=1}^{K}\gamma_{k}^{N,j-1}\psi_{k}(x)+\sum_{k=1}^{b}\gamma
_{k+K}^{N,j-1}\nu_{k}^{N,j-1}(x).\nonumber
\end{align}

\begin{rem}
For definiteness the regression steps (\ref{vp1})-(\ref{gjm1}) are chosen due
to the Tsitsiklis-van Roy (TV)  approach~\cite{J_TV2001}. With a few minor and
obvious changes our reinforced regression approach may be applied to the
Longstaff-Schwartz (LS) method \cite{J_LS2001} as well. Since the details and the
complexity analysis are very similar, we restrict our selves to the TV
approach in this paper.
\end{rem}

\subsection{Spelling out the algorithm}

\label{sec:algorithm2} Let us spell out the above pseudo-algorithm under the
choice (\ref{vb}) of reinforcing functions in more details (general case can
be studied in a similar way). In a pre-computation step we first generate and
save for $m=1,\ldots,N,$ the values%
\begin{equation}
\psi_{k}(Z_{j}^{(m)}),\text{ \ \ }g_{i}(Z_{j}^{(m)}),\text{ \ \ }1\leq j\leq
i\leq\mathcal{J},\text{ \ \ }1\leq k\leq K.\label{pre}%
\end{equation}

\paragraph{Backward procedure}

At the initial time $j=\mathcal{J},$ we set $C_{N,\mathcal{J}}:=0.$ For a
generic backward step $j<\mathcal{J}$ we assume that the quantities%
\begin{equation}
C_{N,j}(Z_{l}^{(m)}),\text{ \ \ }0\leq l\leq j,\text{ \ \ }%
m=1,...,N,\label{jin1}%
\end{equation}
as well as the coefficients $\gamma^{N,j}\in\mathbb{R}^{K+1}$ are already
computed and stored, where formally $C_{N,j}(x)$ satisfies
\begin{align}
C_{N,j}(x)  &  =\sum_{k=1}^{K}\gamma_{k}^{N,j}\psi_{k}(x)+\gamma_{K+1}%
^{N,j}\nu_{1}^{N,j}(x)\label{cf}%
\end{align}
with $\nu_{1}^{N,j} =\max(g_{j+1},C_{N,j+1}).$ Let us now assume that
$0<j\leq\mathcal{J},$ and proceed to time $j-1.$ We first compute (\ref{mpsi})
and (\ref{vp1}). The latter one, $\mathcal{V}_{j},$ is directly obtained by
(\ref{jin1}) for $l=j$ and the pre-computed values (\ref{pre}). To compute
(\ref{mpsi}), we need $\Psi_{K+1}^{N,j-1}(Z_{j-1}^{(m)})=\nu_{1}%
^{N,j-1}(Z_{j-1}^{(m)}),$ $m=1,\ldots, N.$ Hence, we set%
\begin{align*}
\nu_{1}^{N,j-1}(Z_{j-1}^{(m)})  &  =\max(g_{j}(Z_{j-1}^{(m)}),C_{N,j}%
(Z_{j-1}^{(m)}))
\end{align*}
for $m=1,\ldots,N,$ using (\ref{jin1}) for $l=j-1.$ Next we may compute (and
store) the coefficients vector (\ref{gjm1}), i.e., $\gamma^{N,j-1},$ using
(\ref{mpsi}) and (\ref{vp1}), and formally establish (\ref{cf}). In order to
complete the generic backward step, we now need to evaluate
\begin{align}
C_{N,j-1}(Z_{l}^{(m)})=\sum_{k=1}^{K}\gamma_{k}^{N,j-1}\psi_{k}(Z_{l}%
^{(m)})\label{t1}\\
+\gamma_{K+1}^{N,j-1}\nu_{1}^{N,j-1}(Z_{l}^{(m)}),\label{t2}%
\end{align}
for $m=1,...,N,$ $0\leq l\leq j-1.$ The first part (\ref{t1}) is directly
obtained from the pre-computation (\ref{pre}) and the coefficients
(\ref{gjm1}) computed in this step. For the second part (\ref{t2}), we have
that
\begin{align*}
\nu_{1}^{N,j-1}(Z_{l}^{(m)})  &  =\max(g_{j}(Z_{l}^{(m)}),C_{N,j}(Z_{l}%
^{(m)}))
\end{align*}
for $m=1,\ldots,N,$ and $0\leq l\leq j-1.$ Thus the terms (\ref{t2}) are
directly obtained from (\ref{pre}) the coefficients (\ref{gjm1}), and
(\ref{jin1}).

\begin{rem}
(i) Keeping track of the whole set (\ref{jin1}) (rather than some subset, for
example $j-1$ $\leq$ $l$ $\leq$ $j$) in the above procedure is subtle and
necessary due to the nested structure of the additional basis functions
backwardly generated. From a more formal programming point of view this pops
up as a natural ingredient for the logical recursion invariant when going from
$j$ to $j-1.$

(ii) As can be seen, each approximation $C_{N,j-1}$ nonlinearly depends on all
previously estimated continuation functions $C_{N,j},\ldots, C_{N,\mathcal{J}%
-1}$ and hence on all ``features'' $(g_{l}(Z_{l}^{(m)}),\psi_{k}(Z_{l}%
^{(m)}),\, k=1,\ldots, K,\, m=1,\ldots, N,\, l=j,j+1,\ldots, \mathcal{J}).$ In
this sense our procedure finds a sparse nonlinear type approximation for the
continuation functions based on simulated ``features''. Compared to other
nonlinear learning type algorithms (see, e.g., \cite{becker2018deep}), our
procedure doesn't require any nonlinear optimization over high-dimensional
parameter spaces.
\end{rem}

\subsubsection*{Cost estimation}

The total cost needed to perform the pre-computation (\ref{pre}) is about
$\frac{1}{2}N\mathcal{J}^{2}c_{f}+N\mathcal{J}Kc_{f},$ where $c_{f}$ denotes
the maximal cost of evaluating each function $g_{j},$ $j=0,\ldots,
\mathcal{J}, $ and $\psi_{k},$ $k=1,\ldots, K,$ at a given point. The cost of
one backward step from $j$ to $j-1$ can be then estimated from above by
\begin{align*}
&  NK^{2}c_{\ast}\text{ \ \ due to computation of (\ref{gjm1})}\\
&  NKjc_{\ast}\text{ \ \ due to the construction of (\ref{t1})+(\ref{t2}),}%
\end{align*}
where $c_{\ast}$ denotes the sum of costs due to the addition and
multiplication of two reals. Hence the total cost of the above algorithm can
be upper bounded by
\begin{equation}
\frac{1}{2}N\mathcal{J}^{2}c_{f}+N\mathcal{J}Kc_{f}+N\mathcal{J}K^{2}c_{\ast
}+\frac{1}{2}N\mathcal{J}^{2}Kc_{\ast}\label{cost}%
\end{equation}
including the pre-computation.

\subsection{Lower estimate based on a new realization}

Suppose that the backward algorithm of Section~\ref{psalg} has been carried
out, and that we now have an independent set of realizations $(\widetilde
{Z}_{j}^{(m)},$ $j=0,\ldots,\mathcal{J})$ with $\widetilde{Z}_{0}^{(m)}%
=X_{0},$ $m=1,\ldots,N_{\text{test}}.$ In view of (\ref{form1}) and
(\ref{form2}), let us introduce the stopping rule
\begin{equation}
\tau_{N}=\min\bigl\{ j: 1\leq j\leq\mathcal{J},\text{ \ \ }g_{j}(Z_{j})\geq
C_{N,j}(Z_{j})\bigr\}.\label{stop1}%
\end{equation}
A lower estimate of $V_{0}$ is then obtained via%
\begin{equation}
\underline{V_{0}}:=\frac{1}{N_{\text{test}}}\sum_{m=1}^{N_{\text{test}}%
}g_{\tau_{N}^{(m)}}\Bigl(\widetilde{Z}_{\tau_{N}^{(m)}}^{(m)}\Bigr).\label{V0}%
\end{equation}
Here the index $N$ in the $C_{N,j}$ indicates that these objects are
constructed using the simulation sample used in (\ref{psalg}). As a result,
(\ref{stop1}) is a suboptimal stopping time and (\ref{V0}) is a lower biased
estimate. Let us consider the computation of (\ref{stop1}). The coefficient
vectors $\gamma^{N,j},$ $1\leq j\leq\mathcal{J},$ were already computed in the
backward algorithm above. We now have to consider the computation of
$C_{N,j}(Z)$ for an arbitrary point $Z\in\{\widetilde{Z}_{j}^{(m)},\,
m=1,\ldots,N_{\text{test}}\}$ at a particular time $j,$ for $1\leq
j\leq\mathcal{J}.$ For this we propose the following backward procedure.

\subsubsection*{Procedure for computing $C_{N,j}(Z)$ for arbitrary state $Z$}

\begin{enumerate}
\item We first (pre-)compute $\psi_{k}(Z)$ for $1\leq k\leq K,$ and $g_{l}(Z)
$ for $j<l\leq\mathcal{J},$ leading to the cost of order $\left(
K+(\mathcal{J}-j)\right)  c_{f}.$

\item Next compute $C_{N,j}(Z)$ recursively as follows:

\begin{enumerate}
\item Initialize $C_{N,\mathcal{J}}(Z):=0.$ Once $C_{N,l}(Z)$ with
$j<l\leq\mathcal{J},$ is computed and saved, evaluate $\nu_{1}^{N,l-1}(Z)$
using (\ref{vb}).

\item Compute%
\[
C_{N,l-1}(Z)=\sum_{k=1}^{K}\gamma_{k}^{N,l-1}\psi_{k}(Z)+\gamma_{K+1}%
^{N,l-1}\nu_{1}^{N,l-1}(Z)
\]
at a cost of order $Kc_{\ast}.$ In this way we proceed all the way down to
$C_{N,j}(Z),$ at a total cost of $\left(  K+(\mathcal{J}-j)\right)
c_{f}+K\left(  \mathcal{J}-j\right)  c_{\ast}$ including the pre-computation step.
\end{enumerate}
\end{enumerate}

Due to the procedure described above, the costs of evaluating (\ref{V0}),
based on the worst case costs of computing (\ref{stop1}), will be of order%
\begin{equation}
N_{\text{test}}\mathcal{J}Kc_{f}+\frac{1}{2}\mathcal{J}^{2}N_{\text{test}%
}c_{f}+\frac{1}{2}N_{\text{test}}K\mathcal{J}^{2}c_{\ast}.\label{costl}%
\end{equation}
Obviously, (for $N_{\text{test}}=N$) this is the same order as for the
regression based backward induction procedure described in Section~\ref{psalg}.

\subsection{Cost comparison standard vs reinforced regression}
\label{sec:cost-comp}
From the cost analysis of the reinforced regression algorithm it is obviously
inferable that the standard regression procedure, that is, the regression
procedure due to a fixed basis $\psi_{1},\ldots,\psi_{K}$ without
reinforcement, would require a computational cost of order
\begin{equation}
N\mathcal{J}Kc_{f}+N\mathcal{J}K^{2}c_{\ast}\label{costst}%
\end{equation}
for computing the regression coefficients. As an ultimate goal of the
reinforcement method we will try to achieve an accuracy comparable with
standard regression, while the cardinality of the fixed basis is vastly
reduced. If we denote the cardinality of the fixed basis in the reinforced
regression by $K_{r},$ the cost ratio with respect to standard regression is
then given by (\ref{cost})/(\ref{costst}), that is%
\begin{multline*}
\frac{\text{Cost of coefficients the reinforced regression}}{\text{Cost for
coefficients of the standard regression}}\\
=\frac{K_{r}+\mathcal{J}/2}{K}\frac{1+K_{r}c_{\ast}/c_{f}}{1+Kc_{\ast}/c_{f}}.
\end{multline*}
On the other hand, a subsequent lower estimate based on a new realization in
the standard case would require about $N_{\text{test}}\mathcal{J}Kc_{f},$
yielding a cost ratio (see (\ref{costl})),%
\begin{multline*}
\frac{\text{Cost new simulation reinforced regression}}{\text{Cost new
simulation standard regression}}\\
=\frac{K_{r}+\mathcal{J}/2}{K}+\frac{1}{2}\frac{\mathcal{J}K_{r}}{K}c_{\ast
}/c_{f}.
\end{multline*}
From this we conclude that the cost reduction due to the reinforced regression
algorithm is ``large'' when $\left(  K_{r}+\mathcal{J}/2\right)  /K$\ is
``small'', while there is also a ``large'' reduction in the lower bound
construction when in addition $\mathcal{J}c_{\ast}\lesssim$ $c_{f}$ (for example).

\section{Some theoretical results}

Let us consider for a random vector $(X,Y)\in\mathbb{R}^{d}\times\mathbb{R}$
on some probability space $(\Omega,\mathcal{F},\mathbb{P}),$ a problem of
estimating the conditional expectation
\begin{equation}
u(x)=\mathrm{E}\left[  Y\,|\,X=x\right]  ,\label{form}%
\end{equation}
based on a sample $(X^{(n)},Y^{(n)}),$ $n=1,\ldots,N,$ from the joint
distribution of $(X,Y).$ Suppose that the regression basis consists of a fixed
set of standard basis functions $\psi_{k}:$ $\mathbb{R}^{d}\rightarrow
\mathbb{R},$ $k=1,...,K,$ (for example, polynomials) and a set of auxiliary
basis functions $\nu_{1},\ldots,\nu_{b},$ where typically $b$ is much smaller
than $K.$ The idea is that the function $u$ can be well approximated by
functions from $\mathcal{V}_{b}:=\mathsf{span}\left\{  \nu_{1},\ldots,\nu
_{b}\right\} .$ In this case one can consider the least squares problem,%
\begin{equation}
\widetilde{\beta}:=\underset{\beta\in\mathbb{R}^{K+b}}{\arg\inf}\sum_{n=1}%
^{N}\left(  Y^{(n)}-\sum_{k=1}^{K}\widetilde{\beta}_{k}\psi_{k}( X^{(n)})
-\sum_{k=1}^{b}\widetilde{\beta}_{K+k}\nu_{k}( X^{(n)}) \right)
^{2}\label{lsq1}%
\end{equation}
and set
\begin{equation}
\widetilde{u}\left(  x\right)  =\sum_{k=1}^{K}\widetilde{\beta}_{k}\psi
_{k}\left(  x\right)  +\sum_{k=1}^{b}\widetilde{\beta}_{K+k}\nu_{k}\left(
x\right)  .\label{ste}%
\end{equation}
The following theorem provides error bounds for $\widetilde{u},$ see
\cite{audibert2011robust}.

\begin{thm}
\label{rgth} (Accuracy standard global regression) Fix some $\varepsilon
\in(0,1).$ Suppose that
\[
\sup_{x\in\mathbb{R}^{d}}\left\vert u(x)\right\vert \leq L\text{ \ \ and
\ \ }\sup_{x\in\mathbb{R}^{d}}\operatorname{Var}\left[  Y\,|\,X=x\right]
\leq\sigma^{2},
\]
then it holds with probability at least $1-\varepsilon$
\begin{align}
\label{Gyth}\int\left\vert \widetilde{u}(x)-u(x)\right\vert ^{2}\mu(dx)  &
\lesssim\max\left(  \sigma^{2},L^{2}\right)  \frac{\left(  1+\ln N\right)
K+\log(\varepsilon^{-1})}{N}\nonumber\\
&  +\inf_{w\,\in\Psi_{K}+\mathcal{V}_{b}}\int_{\mathbb{R}^{d}}\left\vert
w(x)-u(x)\right\vert ^{2}\mu(dx)
\end{align}
where $\Psi_{K}:=\mathsf{span}\left\{  \psi_{1},\ldots,\psi_{K}\right\} ,$
$\mu$ denotes the distribution of $X$ in (\ref{form}) and $\lesssim$ stands
for inequality up to some absolute constant.
\end{thm}

In view of (\ref{lsq1}) one trivially has for any arbitrary but fixed
$w(x)\in\Psi_{K}+\mathcal{V}_{b},$%
\begin{multline*}
\inf_{\widetilde{\beta}\in\mathbb{R}^{K+b}}\sum_{n=1}^{N}\left(  Y^{(n)}%
-\sum_{k=1}^{K}\widetilde{\beta}_{k}\psi_{k}\left(  X^{(n)}\right)
-\sum_{k=1}^{b}\widetilde{\beta}_{K+k}\nu_{k}\left(  X^{(n)}\right)  \right)
^{2}\\
=\inf_{\widehat{\beta}\in\mathbb{R}^{K+b}}\sum_{n=1}^{N}\left(  Y^{(n)}%
-w(X^{(n)})-\sum_{k=1}^{K}\widehat{\beta}_{k}\psi_{k}\left(  X^{(n)}\right)
\right. \\
\left.  -\sum_{k=1}^{b}\widehat{\beta}_{K+k}\nu_{k}\left(  X^{(n)}\right)
\right)  ^{2}%
\end{multline*}
with the corresponding estimator
\begin{equation}
\widehat{u}( x) =\sum_{k=1}^{K}\widehat{\beta}_{k}\psi_{k}(x) -\sum_{k=1}%
^{b}\widehat{\beta}_{K+k}\nu_{k}(x)\label{up1}%
\end{equation}
of the function $u(x)-w(x).$ Due to (\ref{Gyth}) we thus have for (\ref{up1}),%
\begin{multline*}
\int\left\vert \widehat{u}(x)-u(x)+w(x) \right\vert ^{2}\mu(dx) \lesssim
\max\left(  \sigma^{2},L_{w}^{2}\right)  \frac{\left(  1+\ln N\right)
K+\log(\varepsilon^{-1})}{N}+\delta_{K}%
\end{multline*}
with%
\begin{align*}
L_{w} :=\sup_{x\in\mathbb{R}^{d}}\left\vert u(x)-w(x) \right\vert ,
\quad\delta_{K} :=\inf_{w\,\in\Psi_{K}+\mathcal{V}_{b}}\int_{\mathbb{R}^{d}%
}\left\vert w(x)-u(x)\right\vert ^{2}\mu(dx).
\end{align*}
Since the choice of $w$ was arbitrary, we derive with probability at least
$1-\varepsilon$
\begin{align}
\label{g2}\int\left\vert \widetilde{u}\left(  x\right)  -u(x)\right\vert
^{2}\mu(dx) \lesssim\max\left(  \sigma^{2},L_{\star}^{2}\right)  \frac{\left(
1+\ln N\right)  K+\log(\varepsilon^{-1})}{N}+\delta_{K},\nonumber
\end{align}
where
\[
L_{\star}:=\inf_{w\,\in\Psi_{K}+\mathcal{V}_{b}}\sup_{x\in\mathbb{R}^{d}%
}\left\vert u(x)-w\left(  x\right)  \right\vert .
\]
The reduction of the bound $L$ in (\ref{Gyth}) to $L_{\star}$ is of prime
importance in the backward algorithm developed in Section \ref{sec:algorithm2}%
. In particular, for a diffusion process $X$, the conditional variance of the
underlying process $Z_{j}=X_{t_{j}}$ at $t_{j},$ given its state
$Z_{j-1}=X_{t_{j-1}},$ is of order $O(t_{j}-t_{j-1}).$ It is not difficult to
show that, under some conditions,
\[
\mathsf{Var}[V_{j}(Z_{j})|Z_{j-1}=z]\leq\mathsf{E}\bigl[(V_{j}(Z_{j}%
)-V_{j}(Z_{j-1}))^{2}|Z_{j-1}=z\bigr]=O(t_{j}-t_{j-1}),
\]
uniformly in $z,$ implying $\sigma^{2}\lesssim\max_{j}(t_{j}-t_{j-1})$ in
(\ref{Gyth}). As a result,
\begin{align*}
L_{\star} & \leq \max_{j}\sup_{z}\mathsf{E}[|V_{j}(Z_{j})-V_{j}(z)||Z_{j-1}%
=z]|\\
&  \leq \max_{j}\sup_{z}\sqrt{\mathsf{E}\bigl[|V_{j}(Z_{j})-V_{j}%
(z)|^{2}|Z_{j-1}=z\bigr]}\lesssim\max_{j}\sqrt{t_{j}-t_{j-1}}.
\end{align*}
So in this case $\sigma^{2}\ll L$ in (\ref{Gyth}) and the decrease of $L$ to
$L_{\star}\asymp\sigma$ will result in a substantial computational gain.

\section{Numerical examples}

In this section we illustrate the performance of reinforced regression based
Monte Carlo algorithms by considering two option pricing problems in finance.
\label{num}

\subsection{Bermudan max-call on \(d\) assets}

This is a benchmark example studied in \cite{BG} 
among others. Specifically, the model with $d$ identically distributed
assets is considered, where each underlying has dividend yield $\delta $.
The risk-neutral dynamic of assets is given by
\begin{equation*}
\frac{dX_{t}^{k}}{X_{t}^{k}}=(r-\delta )dt+\sigma dW_{t}^{k},\quad k=1,...,d,
\end{equation*}%
where $W_{t}^{1},\ldots,W_{t}^d$ are independent one-dimensional Brownian
motions and $r,\delta ,\sigma $ are constants. At any time $t\in
\{t_{0},\ldots,t_{\mathcal{J}}\}$ the holder of the option may exercise it and
receive the payoff
\begin{equation*}
g(X_{t})=(\max (X_{t}^{1},\ldots,X_{t}^{d})-K)^{+}.
\end{equation*}%
We take $t_{i}=iT/\mathcal{J},\,i=0,...,\mathcal{J}$, with $T=3,\,\mathcal{J}%
=9$ and $X_{0}=(X_{0}^{1},\ldots ,X_{0}^{d})^{T}$ with $X_{0}^{1}=\ldots=X_{0}^{d}=x_{0}.$   The lower bounds for the standard least-squares approach and the reinforced  regression algorithm are presented in Table~\ref{Tab1} depending on  dimension \(d\) and the choice of basis functions.
In both cases we generate $N=1,000,000$ paths to estimate regression coefficients  and another $ N_{\mathrm{test}}=1,000,000$  to construct  lower bounds (see \eqref{V0}) presented in Table~\ref{Tab1}. The dual upper bounds in the last column of Table~\ref{Tab1} are obtained based on the reinforced regression  using \(1000\) inner paths and \(1,000,000\) outer paths.  
As one can see, there is a clear improvement  in bounds when using the same basis functions across all dimensions.  This improvement is especially pronounced in small dimensions. These results also show  that the RLS algorithm is more efficient than the LS algorithm. Indeed, as can be seen from Table~\ref{Tab1},  the lower bounds for the RLS algorithm achieved when using linear polynomials,  can be obtained  for the LS algorithm only on quadratic ones resulting in a cost reduction of order  \(\frac{2d+\mathcal{J}}{d(d+1)},\) see Section~\ref{sec:cost-comp}.
\begin{table}[h]
\begin{center}
\begin{tabular}
[c]{|c|l|c|c|c|}\hline
\multirow{2}{*}{Dimension} & \multirow{2}{*}{Basis functions} &
\multicolumn{2}{c|}{Lower bounds} & \multirow{2}{*}{Upper bounds} \\ \cline{3-4 }
&  & Regression & Reinf. Reggression &\\\hline
\multirow{3}{*}{2} & $1, X_{i}$ & 12.91(0.018) &13.77(0.015)&14.12(0.042)
\\
& $1, X_{i}, X_{i}X_{j}$ & 13.75(0.014) & 13.86(0.016)&13.97(0.026)
\\
& $1, X_{i}, g(X)$ & 13.66(0.023) & - &14.09(0.071)
\\\hline
\multirow{3}{*}{5} & $1, X_{i}$ &25.25(0.013)  &25.99(0.017)&26.34(0.080)
\\
& $1,  X_{i}, X_{i}X_{j}$ & 25.93(0.020) & 26.12(0.017)&26.22(0.026)
\\
& $1, X_{i}, g(X)$ & 25.82(0.026) & - &26.35(0.081)
\\\hline
\multirow{3}{*}{10} & $1,  X_{i}$ & 37.95(0.025) & 38.22(0.020)&38.48(0.073)
\\
& $1,  X_{i}, X_{i}X_{j}$ & 38.27(0.014) &38.31(0.021)&38.41(0.028)
\\
& $1, X_{i}, g(X)$ & 38.03(0.016)  & - &38.59(0.066)
\\\hline
\multirow{3}{*}{20} & $1,  X_{i}$ & 51.48(0.019)  & 51.61(0.024)&51.88(0.091)
\\
& $1,  X_{i}, X_{i}X_{j}$ & 51.72(0.023) & 51.73(0.023)&51.79(0.035)
\\
& $1, X_{i}, g(X)$ & 51.50(0.020) & - &51.87(0.122)
\\\hline
\end{tabular}
\end{center}\caption{Bounds (with $95\%$ confidence intervals) for the Bermudan max-call
with parameters $K=100,\,r=0.05$, $\protect\sigma =0.2$, $\protect\delta %
=0.1, $ $x_{0}=100$ and different values of $d.$  \label{Tab1}}
\end{table}
The basis \(1, (X_{i}), g(X)\) is skipped for the reinforced regression, since \(g\) is already included (at least at time \(\mathcal{J}-1\)).
\subsection{Bermudan cancelable swap}

We test our algorithm in the case of the so-called complex structured asset
based cancelable swap. We consider a multi-dimensional Black-Scholes model, that is, we define the
dynamic of $d$ assets $X_{l},$ $l=1,\ldots,d,$ under the risk-neutral measure
via a system of SDEs
\[
dX_{l}(t)=(\rho-\delta)X_{l}(t)dt+\sigma_{l}X_{l}(t)dW_{l}(t),\quad0\leq t\leq
T,\quad l=1,\ldots,d.
\]
Here $W_{1}(t),\ldots,W_{d}(t)$ are correlated $d$-dimensional Brownian
motions with time independent correlations $\rho_{lm}=t^{-1}\mathsf{E}%
[W_{l}(t)W_{m}(t)],$ $1\leq l,m\leq d.$ The continuously compounded interest
rate $r$ and a dividend rate $\delta$ are assumed to be constant.
Define the asset based cancelable coupon swap. Let $t_{1},\ldots
,t_{\mathcal{J}}$ be a sequence of exercise dates. Fix a quantile $\alpha,$
$0<\alpha<1$, numbers $1\leq n_{1}<n_{2}\leq d$ (we assume $d\geq2$), and
three rates $s_{1},s_{2},s_{3}$. Let
\[
N(i)=\#\{l:1\leq l\leq d,\ X_{l}(t_{i})\leq(1-\alpha)X_{l}(0)\},
\]
that is, $N(i)$ is the number of assets which at time $t_{i}$ are below
$1-\alpha$ percents of the initial value. We then introduce the random rate
\[
a(i)=s_{1}1_{\left\{  N(i)\leq n_{1}\right\}  }+s_{2}1_{\left\{
n_{1}<N(i)\leq n_{2}\right\}  }+s_{3}1_{\left\{  n_{2}<N(i)\right\}  }%
\]
and specify the $t_{i}$-coupon to be
\[
C(i)=a(i)(t_{i}-t_{i-1}).
\]
For pricing this structured product, we need to compare the coupons $C(i)$
with risk free coupons over the period $[t_{i-1},t_{i}]$ and thus to consider
the discounted net coupon process
\[
\mathcal{C}(i)=e^{-rt_{i}}(e^{r(t_{i}-t_{i-1})}-1-C(i)),\quad i=1,\ldots,
\mathcal{J}.
\]
The product value at time zero may then be represented as the solution of an
optimal stopping problem with respect to the adapted discounted cash-flow,
obtained as the aggregated net coupon process,%
\[
V_{0}=\sup\limits_{\tau\in\{1,\ldots,\mathcal{J}\}}\mathsf{E}[\mathcal{Z}%
_{\tau}], \quad\mathcal{Z}_{j}:=\sum\limits_{i=1}^{j}\mathcal{C}(i).
\]
For our experiments, we choose a five-year option with semiannual exercise
possibility, that is, we have
\[
\mathcal{J}=10,\text{ \ \ }t_{i}-t_{i-1}=0.5,\text{ \ \ }1\leq i\leq10,
\]
on a basket of $d=20$ assets. In detail, we take the following values for the
parameters,
\[%
\begin{split}
d  &  =20,\quad r=0.05,\quad\delta=0,\quad\sigma_{l}=0.2,\quad X_{l}%
(0)=100,\quad1\leq l,m\leq20,\\
d_{1}  &  =5,\quad d_{2}=10,\quad\alpha=0.05,\quad s_{1}=0.09,\quad
s_{2}=0.03,\quad s_{3}=0,
\end{split}
\]
and
\[
\rho_{lm}=%
\begin{cases}
\rho, & l\neq m,\\
1, & l=m.
\end{cases}
\]
As to the basis functions, we used a constant, the discounted net coupon
process $\mathcal{C}(i)$ and the order statistics $X_{(1)}\leq X_{(2)}%
\leq\ldots\leq X_{(n)}$. Table~\ref{my-label} shows the results of the
numerical experiment comparing the lower and the corresponding dual upper
bounds by the standard linear regression method with fixed basis (the second
column of Table~\ref{my-label}) and by the reinforced regression approach
described in Section~\ref{sec:algorithm2} with one additional basis function
$(\nu_{1}^{N,j}).$ The main conclusion is that the reinforced regression
algorithm delivers estimates of the same quality as the standard least squares
approach by using much less basis functions (sparse basis). As a result the
new algorithm turns out to be computationally cheaper. \begin{table}[h]
{\scriptsize \label{my-label} }
\par
\begin{center}%
\begin{tabular}
[c]{|c|l|c|c|}\hline
\multirow{2}{*}{\(\rho\)} & \multirow{2}{*}{Basis functions} &
\multicolumn{2}{c|}{Regression}\\\cline{3-4 }
&  & Low Estimation & High Estimation\\\hline
\multirow{2}{*}{0} & $1, \mathcal{C}, X_{(i)}$ & 171.59(0.037) &
177.24(0.061)\\
& $1, \mathcal{C}, X_{(i)}, X_{(i)}X_{(j)}$ & 173.62(0.044) &
177.33(0.062)\\\hline
\multirow{2}{*}{0.2} & $1, \mathcal{C}, X_{(i)}$ & 180.0(0.060) &
199.62(0.125)\\
& $1, \mathcal{C}, X_{(i)}, X_{(i)}X_{(j)}$ & 188.01(0.055) &
197.02(0.143)\\\hline
\multirow{2}{*}{0.5} & $1, \mathcal{C}, X_{(i)}$ & 176.43(0.073) &
201.21(0.189)\\
& $1, \mathcal{C}, X_{(i)}, X_{(i)}X_{(j)}$ & 183.41(0.033) &
196.58(0.147)\\\hline
\multirow{2}{*}{0.8} & $1, \mathcal{C}, X_{(i)}$ & 133.29(0.065) &
158.12(0.197)\\
& $1, \mathcal{C}, X_{(i)}, X_{(i)}X_{(j)}$ & 140.17(0.061) &
153.49(0.106)\\\hline
\end{tabular}
\par%
\begin{tabular}
[c]{|c|l|c|c|}\hline
\multirow{2}{*}{\(\rho\)} & \multirow{2}{*}{Basis functions} &
\multicolumn{2}{c|}{Reinf. regression}\\\cline{3-4 }
&  & Low Estimation & High Estimation\\\hline
\multirow{2}{*}{0} & $1, \mathcal{C}, X_{(i)}$ & 173.28(0.031) &
177.32(0.091)\\
& $1, \mathcal{C}, X_{(i)}, X_{(i)}X_{(j)}$ & 174.33(0.036) &
176.58(0.057)\\\hline
\multirow{2}{*}{0.2} & $1, \mathcal{C}, X_{(i)}$ & 187.57(0.057) &
195.09(0.121)\\
& $1, \mathcal{C}, X_{(i)}, X_{(i)}X_{(j)}$ & 188.07(0.046) &
195.95(0.108)\\\hline
\multirow{2}{*}{0.5} & $1, \mathcal{C}, X_{(i)}$ & 181.98(0.047) &
194.04(0.088)\\
& $1, \mathcal{C}, X_{(i)}, X_{(i)}X_{(j)}$ & 183.93(0.057) &
194.97(0.127)\\\hline
\multirow{2}{*}{0.8} & $1, \mathcal{C}, X_{(i)}$ & 138.41(0.087) &
153.08(0.106)\\
& $1, \mathcal{C}, X_{(i)}, X_{(i)}X_{(j)}$ & 139.62(0.035) &
152.57(0.096)\\\hline
\end{tabular}
\end{center}
\caption{Comparison of the standard linear regression method and the
reinforced regression algorithm for the problem of pricing cancelable swaps}%
\end{table}
\pagebreak
\section*{References}

\bibliographystyle{plain}
\bibliography{deep_optimal_stop-1}

\end{document}